\documentclass[12pt]{article}

\usepackage{amstext    }
\usepackage{amsthm    }
\usepackage{a4}
\usepackage[mathscr]{eucal} 
\usepackage{mathrsfs}

\usepackage{amsmath}
\usepackage{amssymb}
\usepackage{amscd}

\numberwithin{equation}{section}

\newtheorem{theorem}{Theorem}[section]

\newtheorem{proposition}[theorem]{Proposition}
\newtheorem{corollary}[theorem]{Corollary}
\newtheorem{lemma}[theorem]{Lemma}
\newtheorem{remark}[theorem]{Remark}

\newcommand{\cali}[1]{\mathscr{#1}} 

\newcommand{\FS}{{\rm FS}}

\newcommand{\volume}{{\rm volume}}
\newcommand{\lov}{{\rm lov}}
\newcommand{\supp}{{\rm supp}}
\newcommand{\const}{{\rm const}}
\newcommand{\dist}{{\rm distance}}

\newcommand{\loc}{{loc}}
\newcommand{\ddc}{{\rm dd}^{\rm c}}
\newcommand{\dc}{{\rm d}^{\rm c}}
\renewcommand{\d}{{\rm d}}
\newcommand{\dbar}{\overline\partial}

\newcommand{\id}{{\rm id}}

\newcommand{\Jac}{{\rm Jac}}
\newcommand{\Aut}{{\rm Aut}}
\newcommand{\Leb}{{\rm Leb}}

\newcommand{\Ac}{\cali{A}}

\newcommand{\Cc}{\cali{C}}
\newcommand{\Dc}{\cali{D}}

\newcommand{\Lc}{\cali{L}}

\newcommand{\Rc}{\cali{R}}

\newcommand{\C}{\mathbb{C}}

\renewcommand{\P}{\mathbb{P}}

\title{Attracting current and equilibrium measure for attractors
  on $\P^k$}
\author{Tien-Cuong Dinh}

\begin{document}

\maketitle

\begin{abstract} Let $f$ be a holomorphic
  endomorphism of $\P^k$ having an attracting set $\Ac$. 
We construct an attracting current and an equilibrium
measure associated to $\Ac$. The attracting 
current is weakly laminar and extremal in the cone of invariant currents. 
The equilibrium measure is mixing and has maximal entropy.
\end{abstract}

\section{Introduction}

Let $f$ be a holomorphic endomorphism of $\P^k$ of algebraic degree
$d\geq 2$. 
If $S$ is a smooth positive closed $(1,1)$-form of mass 1 then
$d^{-n}(f^n)^*S$ converges weakly to a positive closed $(1,1)$-current $T$ of
  mass 1 which does not depend on $S$ and has locally H{\"o}lder 
  continuous potentials, see Forn\ae ss-Sibony
  \cite{FornaessSibony1,Sibony}.  
The current $T$ is totally
  invariant :  $d^{-1}f^*(T)=d^{-k+1}f_*(T)=T$. This is {\it the
    Green current associated to $f$};
its support is
  {\it the Julia set}.
The 
self-intersection $T^m$ of $T$, $1\leq m\leq k$, is {\it the Green current of bidegree $(m,m)$}. It
satisfies $d^{-m}f^*(T^m)=d^{-k+m}f_*(T^m)=T^m$. We refer to
\cite{Federer, Demailly, HarveyShiffman,
  FornaessSibony3, DinhSibony5}
for the basics on the theory of currents and their intersections.
We also obtain $T^m$ as the limit in the sense of current of $d^{-mn}(f^n)^*[I]$,  
where $[I]$ is the current of integration
on a generic projective subspace $I$ of dimension $k-m$
\cite{FornaessSibony4, RussakovskiiShiffman, DinhSibony6}.

Following Yomdin \cite{Yomdin} and Gromov \cite{Gromov}, $f$ has
topological entropy $k\log d$ and the variational principle \cite{KatokHasselblatt,Walters}
implies that all measures, invariant by $f$, have entropy $\leq k\log
d$. 
The
equilibrium measure $\mu:=T^k$ of $f$ is the  unique invariant probability measure of
maximal entropy $k\log d$, see Forn\ae ss-Sibony
\cite{FornaessSibony4} and Briend-Duval \cite{BriendDuval}.
The main part of the dynamics is, in some sense,
concentrated on the support of $\mu$. In this paper, we study the dynamics
on some attracting sets of $f$. They are disjoint from 
$\supp(\mu)$ or disjoint from $\supp(T^m)$ for some
$m$.

A compact subset $\Ac$ is an {\it attracting set} for $f$
if it admits a neighbourhood $U$ such that $f(U)\Subset U$ and
$\Ac=\cap_{n\geq 0} f^n(U)$. Attracting sets are invariant : $f(\Ac)=\Ac$. If
moreover $f$ admits a dense orbit in $\Ac$ then $\Ac$ is said to be
{\it an attractor}.  We refer to \cite{FornaessWeickert,
  JonssonWeickert, FornaessSibony2, Rong} for basic properties and
examples of attracting sets and attractors, and to \cite{deThelin2} for
the study of weakly saddle measures outside $\supp(\mu)$. 

We now describe the class of attracting sets that we study here.
This class contains the examples studied in the previous references
\footnote{For attracting sets and
attractors in a weaker sense, see \cite{BonifantDabijaMilnor} and the
references therein}. 
Let $I$ and $L$ be two projective subspaces of $\P^k$,  of dimensions $p-1$
and $k-p$ respectively, $1\leq p\leq k-1$. We assume that $I\cap
U=\varnothing$ and that $L$
is contained in $U$. In particular, $L$ does not intersect $I$. 
Note that in our study the case $p=k$ is not interesting. Indeed, $f$ contracts the
Kobayashi pseudo-distance, which is a distance in this case, hence
$\Ac$ is reduced to a point, see also \cite{Sibony}. 

Let $\pi:\P^k\setminus I\rightarrow L$ be the canonical projection of
center $I$. More precisely, if $I(x)$ is the projective space, of dimension $p$, containing $I$
and passing through a point $x\in \P^k\setminus I$ then
$\pi(x)$ is the unique intersection point of $L$
with $I(x)$.  We consider the point
$\pi(x)$ as the origin of the complex vector space $I(x)\setminus I\simeq \C^p$, where 
$I$ is viewed as the hyperplane at infinity of $I(x)\simeq \P^p$.
In other words, $\P^k\setminus I$ is considered as a vector bundle
over $L$. If $x\in L$ we have $\pi(x)=x$.
Our main hypothesis is

\begin{center}
 {\it the open set  $U\cap I(x)$ in
  $I(x)\setminus I\simeq \C^p$ is  
star-shaped at $x$ for every $x\in L$}. 
\end{center}
\noindent
In particular, $U$ and $\P^k\setminus U$ are connected. 
Observe that the previous geometric condition is often easy to check and is stable
under small pertubations on $f$. If $I$ is a point and $L$ is a
projective hyperplane (i.e. $p=1$), the previous
hypothesis is equivalent to the property that the open subset
$\P^k\setminus U$ of $\P^k\setminus L\simeq \C^k$ is
star-shaped at $I$. In particular, 
this last property is satisfied when $\P^k\setminus U$ is convex.

Consider a generic subspace $I'$ of dimension $p-1$
close enough to $I$. Then $d^{-(k-p+1)n}(f^n)^*[I']$ converges to
$T^{k-p+1}$. Since  $d^{-(k-p+1)n}(f^n)^*[I']$ is supported on
$\P^k\setminus U$, the support of $T^{k-p+1}$ does not intersect
$U$. In particular, we have $\supp(T^{k-p+1})\cap \Ac=\varnothing$,
hence 
$\supp(\mu)\cap \Ac=\varnothing$. Studying the dynamics on $\Ac$ allows
us to understand the dynamics outside the support of $\mu$. We will
associate to $\Ac$ a dynamically interesting invariant current and  measure.

The positive closed current $[L]$ of
integration on $L$ is of bidegree $(p,p)$ and of bidimension
$(k-p,k-p)$. It has compact support in $U$. Using a standard convolution
we can construct smooth positive closed $(p,p)$-forms with compact
support in $U$, see also Section \ref{section_structure}.  

\begin{theorem} \label{th_attracting_current}
Let $f$, $I$ and $U$ be as above.
Let $R$ be a continuous positive closed $(p,p)$-form of mass $1$ with 
compact support in $U$. Then
$d^{-(k-p)n} (f^n)_*R$ converges weakly to a positive closed current $\tau$
of mass $1$ with support on $\Ac$. 
The current $\tau$ does not depend on $R$.
Moreover, it
is weakly laminar, invariant ($f_*\tau = d^{k-p}\tau$) and  
is extremal in the cone of invariant positive closed $(p,p)$-currents 
with support in $\Ac$.
\end{theorem}

We call $\tau$ {\it the attracting current associated to $\Ac$}.
The hypothesis on the geometry of $U$ is necessary. 
For this, we can consider the map $f:\P^2\rightarrow\P^2$,
$f[z_0:z_1:z_2]=[z_0^d:z_1^d:z_2^d]$ and $\P^2\setminus U :=
\{{1\over 2} |z_i|\leq |z_j|\leq 2|z_i|, \ \mbox{for } 0\leq i,j\leq
2\}$. Then 
$\Ac=\{z_0=0\}\cup\{z_1=0\}\cup\{z_2=0\}$. We can see that $\Ac$
admits three ``attracting'' currents (currents of integration on
$\{z_i=0\}$). This property is stable under small pertubations of $f$.

There are obvious situation where Theorem \ref{th_attracting_current}
applies. Let $f$ be a polynomial map in $\C^k$ of algebraic degree
$d>1$ with holomorphic extension to $\P^k$, then the hyperplane at
infinity is attracting and Theorem \ref{th_attracting_current} can be
applied to small pertubations of $f$.  

Recall here that the projective space $\P^k$ is endowed with the
Fubini-Study form $\omega_\FS$. This is a K{\"a}hler 
form normalized by $\int\omega_\FS^k=1$. {\it The mass} of a
positive closed $(p,p)$-current $R$ on an open set $W$ is defined by 
$\|R\|_W:=\int_W R\wedge \omega_\FS^{k-p}$. The mass $\|R\|:=\|R\|_{\P^k}$
of $R$ (on $\P^k$)  depends only on its cohomology class in
$H^{p,p}(\P^k,\C)\simeq \C$. 
In particular, if $R$ is of mass 1, it is cohomologous to
$\omega_\FS^p$. Then $\|(f^n)^*R\|=\|(f^n)^*\omega_\FS^p\|=d^{pn}$ and
\begin{eqnarray*}
\|(f^n)_*R\| & = & \int(f^n)_*R\wedge \omega_\FS^{k-p}=\int R\wedge
(f^n)^*\omega_\FS^{k-p}\\
& = & \int \omega_\FS^p\wedge (f^n)^*\omega_\FS^{k-p}=\|(f^n)^*\omega_\FS^{k-p}\|
=d^{(k-p)n}.
\end{eqnarray*}

The current $\tau$ is {\it weakly laminar or web-like}, see \cite{Dinh}, if it
can be decomposed into currents of integration on complex manifolds of
dimension $k-p$ as follows  
\begin{enumerate}
\item[{\it i})] There exists a family of 
complex manifolds $V\subset \P^k$ of dimension $k-p$, not necessarily closed in $\P^k$, with finite
volume. We denote by $[V]$ the current of integration on $V$, it is
positive but not necessarily closed.
\item[{\it ii})] There
is a positive measure $\lambda$ defined on the family of these
currents $[V]$ such
that
$$\int \volume(V)\d\lambda(V)<\infty\qquad \mbox{and}\qquad \tau=\int
[V] \d\lambda(V).$$ 
More precisely, if $\Phi$ is a test smooth $(k-p,k-p)$-form then 
$$\langle \tau,\Phi\rangle = \int \left(\int_V\Phi\right)
\d\lambda(V).$$
\end{enumerate}
Theorem \ref{th_attracting_current} implies that $\Ac$
contains pieces of complex
manifolds of dimension $k-p$. This is an important property which
allows  
to study the dynamics on $\Ac$ using geometric methods, see
\cite{BedfordLyubichSmillie, Dujardin, 
deThelin2} and Section \ref{section_remark}. 

Following an idea of Sibony, see e.g. \cite{Sibony} and
\cite{BedfordSmillie} for historical comments,  we 
define {\it  the equilibrium measure associated to $\Ac$}, as
intersection of invariant currents, by 
$$\nu:=T^{k-p}\wedge\tau.$$ 
This is a probability measure invariant by $f$ ($f_*(\nu)=\nu$) and supported on $\Ac$.
Since $T$ has locally H{\"o}lder continuous
potentials, the measure $\nu$ has positive Hausdorff dimension and has
no mass  on points, see
\cite{BriendDuval, Sibony, DinhSibony4}. 

Using an idea of de Th{\'e}lin \cite{deThelin1}, we prove in
Section \ref{section_entropy} that the topological entropy
of $f$ on $U$ is bounded by $(k-p)\log d$. The variational principle
\cite{KatokHasselblatt, Walters} implies that any invariant measure with support in $\Ac$, in
particular $\nu$, has entropy $\leq (k-p)\log d$.
We have the following theorem.

\begin{theorem} \label{th_measure}
Under the above notation, the measure $\nu$ 
has maximal entropy
  $(k-p)\log d$ and is mixing. 
\end{theorem}

Note that we can associate to $\Ac$ the invariant currents $T^m\wedge \tau$, $1\leq
m\leq k-p-1$, which should be useful in the study of the dynamics on
$\Ac$. 

The precise outline of the paper is as follows. In Section
\ref{section_entropy}, we generalize a theorem due to de Th{\'e}lin to the
case of higher dimension. In Section \ref{section_structure}, we
introduce the main tools needed in the construction of the
attracting current. We use here a new method introduced by Nessim
Sibony and the author in \cite{DinhSibony5} to deal with currents of
bidegree $(p,p)$.
Theorems \ref{th_attracting_current} and \ref{th_measure} are proved in Sections
\ref{section_current} and \ref{section_measure}. Further remarks and
questions are given in Section \ref{section_remark}.
Appendix
\ref{appendix_slicing} 
contains some useful results on the
theory of slicing of currents. Throughout the paper, $\Cc(U)$ denotes
the set of positive closed currents of mass 1 with compact support in
$U$ and $\Dc$ denotes the set of all the limit values of $d^{-(k-p)n}
(f^n)_*R_n$ with $R_n\in\Cc(U)$.


\section{Entropy} \label{section_entropy}

In this section we give a sharp upper bound for the topological entropy of
$f$ outside the support of $T^m$. We follow an idea due to de Th{\'e}lin
in \cite{deThelin1} where he studied the case of dimension 2.

Consider a subset $W$ of $\P^k$, not
necessarily invariant. A family $F$ of points in $W$ is said to be
{\it $(n,\epsilon)$-separated},
$\epsilon>0$, if for all distinct points $x$ and $y$ in $F$, there is
an integer 
$j$, $0\leq j\leq n-1$, such that
$\dist(f^j(x),f^j(y))>\epsilon$. Define
the {\it topological entropy} of $f$ on $W$ by
$$h_t(f,W):=\sup_{\epsilon>0}\limsup_{n\rightarrow\infty} {1\over n} \log
\max\big\{\#F,\ F\subset W \ (n,\epsilon)\mbox{-separated}\big\},$$
see \cite{Bowen, deThelin1}.
By Yomdin \cite{Yomdin} and Gromov \cite{Gromov}, the topological entropy
$h_t(f):=h_t(f,\P^k)$ of $f$ is equal to $k\log d$, see also
\cite{DinhSibony2, DinhSibony3}. The main result of this section is
the following theorem, see \cite{deThelin1} for the case $k=m=2$.

\begin{theorem} \label{th_entropy}
Let $f$ be a holomorphic endomorphism of $\P^k$ of algebraic degree
$d\geq 2$. Let $T$ be the Green current of $f$. If $W$ is a subset of
$\P^k$ such that $\overline W\cap \supp(T^m)=\varnothing$, $1\leq m\leq k$, then 
$$h_t(f,W)\leq (m-1)\log d.$$
\end{theorem}

\proof
The proof is based on an idea of Gromov \cite{Gromov} and on the
speed of convergence toward the Green current by Forn\ae ss-Sibony
\cite{FornaessSibony1, Sibony}.
Fix an open neighbourhood $W'$ of $\overline W$ such that  
$W'\Subset \P^k\setminus \supp(T^m)$.

Let $\Gamma_{[n]}$ denote the graph of the
map
$(f,f^2,\ldots,f^{n-1})$ in $(\P^k)^n$, i.e. the set of points
$$(x,f(x),\ldots,f^{n-1}(x)),\quad x\in\P^k.$$
The projection $\Pi$ on the first factor $\P^k$
defines a biholomorphic map between $\Gamma_{[n]}$ and $\P^k$. 
Consider an $(n,\epsilon)$-separated subset $F\subset W$. We
associate to $F$ the family of balls in  $(\P^k)^n$ of
center $(x,f(x),\ldots,f^{n-1}(x))$, with $x\in F$, and of radius $\epsilon/2$.
Here, the distance between two points
in $(\P^k)^n$ is
the maximum of the distances between their projections on the
factors $\P^k$.
Since $F$ is $(n,\epsilon)$-separated, these balls are disjoint.  
They are contained in $\Pi^{-1}(W')$ when $\epsilon$
is small enough.

An inequality of Lelong \cite{Lelong} implies that 
the volume of $\Gamma_{[n]}$ in each ball is larger than $c\epsilon^{2k}$ for some
constant $c>0$. 
We use here the volume with respect to the metric induced by the K{\"a}hler form 
$\sum\Pi_i^*(\omega_\FS)$, where $\Pi_i$ denotes the canonical projection of
$(\P^k)^n$ on its $i$-th factor.
We deduce an estimate of the number of balls and obtain the following Gromov's inequality
$$h_t(f,W)\leq \lov(f,W'):=\limsup_{n\rightarrow\infty} {1\over n}
\log\volume(\Pi^{-1}(W')\cap \Gamma_{[n]}).$$

We will show that $\volume(\Pi^{-1}(W')\cap \Gamma_{[n]})\lesssim
n^kd^{(m-1)n}$ which implies the theorem. Since $\Gamma_{[n]}$ has
dimension $k$ and since $\Pi$ defines a biholomorphic map between
$\Gamma_{[n]}$ and $\P^k$, we have
\begin{eqnarray*}
\volume(\Pi^{-1}(W')\cap \Gamma_{[n]}) & = &
{1\over k!}\int_{\Pi^{-1}(W')\cap \Gamma_{[n]}} \Big(\sum_{i=0}^{n-1} \Pi_i^*(\omega_\FS)\Big)^k
  \\
& = & {1\over k!}\int_{W'} \Big(\sum_{i=0}^{n-1} (f^i)^*(\omega_\FS)\Big)^k \\
& = & {1\over k!}\sum_{0\leq n_i \leq n-1} \int_{W'}
(f^{n_1})^*(\omega_\FS)\wedge \ldots \wedge (f^{n_k})^*(\omega_\FS).
\end{eqnarray*}
Hence, it is enough to show for $0\leq n_i\leq n$ that
$$\int_{W'}(f^{n_1})^*(\omega_\FS)\wedge \ldots \wedge
(f^{n_k})^*(\omega_\FS) \lesssim d^{(m-1)n}.$$
To this end, we prove by induction on $(r,s)$, $0\leq r\leq m$ and $0\leq s\leq 
k-m+r$, that 
$$\|T^{m-r}\wedge  (f^{n_1})^*(\omega_\FS)\wedge\ldots\wedge
(f^{n_s})^*(\omega_\FS)\|_{W_{r,s}}\leq c_{r,s}d^{n(r-1)},$$
where $W_{r,s}$ is a neighbourhood of $\overline W'$
and $c_{r,s}\geq 0$ is a constant independent of $n$. We obtain the
result by taking $r=m$ and $s=k$.

It is clear that the previous inequality holds for $r=0$ or $s=0$ with
$W_{r,s}=\P^k\setminus\supp(T^m)$ and $c_{r,s}=1$. 
Assume this inequality for the cases of $(r-1,s-1)$ and $(r,s-1)$.
Let $W_{r,s}$ be a neighbourhood of $\overline W'$ strictly contained
in $W_{r-1,s-1}$ and $W_{r,s-1}$. Let $\chi\geq 0$ be a smooth cut-off
function with support in $W_{r-1,s-1}\cap W_{r,s-1}$ which is equal to 1
on $W_{r,s}$. We only have to prove that
$$\int T^{m-r}\wedge (f^{n_1})^*(\omega_\FS)\wedge\ldots\wedge
(f^{n_s})^*(\omega_\FS) \wedge \chi\omega_\FS^{k-m+r-s}\leq c_{r,s}d^{n(r-1)}.$$

Since $T$ and $\omega_\FS$ have the same mass, they belong to the same
cohomology class in $H^{1,1}(\P^k,\C)$. Moreover, $T$ has
locally continuous potentials. Hence, there is a continuous
quasi-psh function
$g$ such that $\omega_\FS=T-\ddc g$. 
We have
$$(f^{n_1})^*\omega_\FS=d^{n_1}T-\ddc (g\circ f^{n_1}).$$
The integral that we have to estimate, is equal to the sum of the following integrals 
$$ d^{n_1}\int T^{m-r+1}\wedge (f^{n_2})^*(\omega_\FS)\wedge\ldots\wedge
(f^{n_s})^*(\omega_\FS) \wedge \chi\omega_\FS^{k-m+r-s}$$
and 
$$ -\int T^{m-r}\wedge \ddc(g\circ f^{n_1})\wedge (f^{n_2})^*(\omega_\FS)\wedge\ldots\wedge
(f^{n_s})^*(\omega_\FS) \wedge \chi\omega_\FS^{k-m+r-s}.$$
Using the case of $(r-1,s-1)$ we can bound the first integral by $c
d^{n(r-1)}$. Stokes' theorem implies that second integral is equal to 
$$ -\int T^{m-r}\wedge (f^{n_2})^*(\omega_\FS)\wedge\ldots\wedge
(f^{n_s})^*(\omega_\FS) \wedge (g\circ f^{n_1})\ddc\chi\wedge
\omega_\FS^{k-m+r-s}$$
which is bounded by
$$ \|g\|_\infty\|\chi\|_{\Cc^2}\|T^{m-r}\wedge (f^{n_2})^*(\omega_\FS)\wedge\ldots\wedge
(f^{n_s})^*(\omega_\FS)\|_{W_{r,s-1}}$$
since $\chi$ has support in $W_{r,s-1}$.
We obtain the result using the case of $(r,s-1)$.
\endproof

From the variational principle \cite{KatokHasselblatt, Walters}, we deduce the following
result.

\begin{corollary} \label{cor_entropy}
Let $\nu$ be a probability measure invariant by $f$. Assume that $\nu$
has compact
support in $\P^k\setminus\supp(T^m)$. Then
the entropy of $\nu$ satisfies 
$$h(\nu)\leq (m-1)\log d.$$ 
\end{corollary}

We refer to \cite{KatokHasselblatt, Walters} for the definition of the
entropy of an invariant measure.


\section{Structure of the cone of currents} \label{section_structure}

In this section, we introduce the main tools that we use in order to
construct the attracting current. The approach is similar to the one in \cite{DinhSibony5}.
Let $\Cc(U)$ denote the set of positive closed $(p,p)$-currents $R$ 
of mass 1 in $\P^k$ with compact support in the open set
$U$. 
Consider a real smooth $(k-p,k-p)$-form $\Phi$ 
on $U$, not necessarily with compact support, 
such that $\ddc\Phi\geq 0$. Define the function $\Lambda_\Phi$ on $\Cc (U)$ by
$\Lambda_\Phi(R):=\langle R,\Phi\rangle$. 
These $\Lambda_\Phi$ play
the role of {\it psh functions} on $\Cc (U)$, see Lemma \ref{lemma_structure_coherent} below.
The following lemma implies that  they separate the currents in $\Cc (U)$.
Recall that $\overline U$ does not intersect the space $I$ of
dimension $p-1$.

\begin{lemma} \label{lemma_density}
Let $\Phi$ be a real smooth $(k-p,k-p)$-form with compact support in
$U$. Then there are smooth $(k-p,k-p)$-forms $\Phi^+$, $\Phi^-$, not necessarily
with compact support, such that $\ddc\Phi^\pm\geq 0$ and
$\Phi=\Phi^+-\Phi^-$. Moreover, there exists a real $(k-p,k-p)$-form $\Psi$ on $U$
such that $\ddc\Psi$ is strictly positive.
\end{lemma}
\proof
If $\Psi$ is as above,
then $\Phi^+:=\Phi+c\Psi$ and
$\Phi^-:=c\Psi$, $c>0$ large enough, satisfy the lemma. We only have to construct such a form
$\Psi$. Fix a point $x$ in $\overline U$. It is enough to construct a
form $\Psi_x$ such that $\ddc\Psi_x\geq 0$ in a neighbourhood of
$\overline U$, and $\ddc \Psi_x$ strictly
positive at $x$. We obtain $\Psi$ by taking a finite sum of such forms
$\Psi_x$.

Let $\varphi$ be a smooth strictly psh function on
$I(x)\setminus I\simeq \C^p$, see the Introduction. Then $\varphi[I(x)]$ defines
a $(k-p,k-p)$-current $\Psi_x'$  on $\P^k\setminus I$ 
such that $\ddc\Psi_x'\geq 0$
and $\ddc\Psi_x'$ does not vanish at
 $x$. We use here the standard method to regularize $\Psi_x'$.

Let $\Aut(\P^k)$ be the group of holomorphic automorphisms of $\P^k$. This is a complex Lie group of
dimension $k^2+2k$. Let $\rho$ be a smooth probability measure on
$\Aut(\P^k)$ with support in a small neighbourhood of the identity
map $\id\in\Aut(\P^k)$. 
Define $$\Psi_x'':=\int A_*(\Psi_x')\d\rho(A).$$
Then $\Psi_x''$ is smooth and is defined out of a small neighbourhood
of $I$. We also have
$$\ddc\Psi_x'':=\int A_*(\ddc\Psi_x')\d\rho(A).$$
Hence, $\ddc\Psi_x''$ is positive 
and $x$ is very close to the support of
$\ddc\Psi_x''$. Replacing $\Psi_x''$ by $A^*(\Psi''_x)$, where $A$ is a
suitable automorphism close to the identity, we can assume that
$\ddc\Psi_x''$ does not vanish at $x$. We can choose a finite number of
automorphisms $A_i$, close to the identity, which fix the point $x$, such that if
$$\Psi_x:=\sum A_i^*(\Psi_x'')$$
then $\ddc \Psi_x$ is strictly positive at $x$.
\endproof

Let $V$ be a complex manifold of dimension $m$. Let $\pi_V: V\times
U\rightarrow V$ and $\pi_U:V\times U\rightarrow U$ be the canonical projections.
Consider a positive closed $(p,p)$-current $\Rc$ in $V\times U$ such
that $\pi_U(\supp(\Rc))\Subset U$. We say that $\Rc$ is {\it
  horizontal}, see also \cite{DinhSibony5}. 
It is shown in Proposition \ref{prop_slice} in
Appendix A that the slice $\langle
\Rc,\pi_V,\theta\rangle$ exists for every $\theta\in V$. This is a
positive closed $(p,p)$-current with compact support in
$\{\theta\}\times U$. We often identify it to a $(p,p)$-current $R_\theta$ in $U$.  
The mass of $R_\theta$ does not depend on $\theta$, see Proposition
\ref{prop_slice_mass}. 
We assume that this mass is
equal to 1.
Then one obtains a map
$h: \theta\mapsto R_\theta$ from $V$ into $\Cc(U)$. We say that the
map $h$ or the family  $(R_\theta)$ defines a {\it structural variety} of $\Cc(U)$.
We deduce from Remark \ref{rk_slice} the following lemma.

\begin{lemma} \label{lemma_structure_coherent}
Let 
$(R_\theta)_{\theta\in V}$ be a structural variety as above. 
Let $\Phi$ be a real smooth $(k-p,k-p)$-form on $U$ such that
$\ddc\Phi\geq 0$.
Then
$\theta\mapsto \langle R_\theta,\Phi\rangle$ defines a psh function on
$\theta\in V$. 
\end{lemma}

In what follows we construct and use only some special structural discs
$(R_\theta)$ such that $R_\theta$ 
depends continuously on $\theta$, where $V$ is a holomorphic disc in
$\C$. The construction of these discs is as follows. 
Fix a chart $W$ of $\Aut(\P^k)$ containing $\id\in\Aut(\P^k)$ and local
holomorphic coordinates $y$, $\|y\|<1$,  
such that $y=0$ at $\id$.
Then we can consider the scalar multiplication $\lambda_\theta(y):=(1-\theta) y$
for $\theta\in\C$ and $\|y\|<\min(1,|1-\theta|^{-1})$. 
Fix  a small neighbourhood $W'\Subset W$ of $\id$. Sometimes we
identify $y$ with the automorphism it represents.

We define a holomorphic self-map $A_\theta$ on $\P^k\setminus
I$. If $x$ is a point in   $\P^k\setminus I$, it belongs to 
$I(x)\setminus I$ and we consider $I(x)\setminus I$ as a complex vector space with
origin at $\pi(x)=I(x)\cap L$, see the Introduction. 
We define $A_\theta(x)$ as the multiplication,  in $I(x)\setminus I$,
of $x$ by $\theta$.
If $\theta=0$ then $A_0$ is equal to the
projection $\pi$ of center $I$ onto $L$.
Otherwise $A_\theta$ can be extended to an automorphism on
$\P^k$ which fixes the points in $I$. In particular, we have $A_1=\id$.

Let $U'\Subset U$ be an open set such that $f(U)\Subset U'$. 
Let $V$ be a simply connected neighbourhood of the interval $[0,1]$ in
$\C$. 
We choose $V$ and $W'$ small
enough so that 
$\lambda_\theta(A)\circ A_\theta (U')\Subset U$ for every $A\in\overline W'$ and 
$\theta\in \overline V$. Here we use the geometric assumption on  $U$
(see the Introduction), which implies this property for $\theta\in[0,1]$.

Consider $R\in \Cc (U')$ and an automorphism $A\in W'$. 
In order to get the picture, the reader can consider the case where
$R$ and the currents 
$R_{A,\theta}$, $\Rc_A$ below, are defined by integration on complex
manifolds.
The currents $R_{A,\theta}$, $\theta\in V$, are obtained as images of
$R$ under a holomorphic family of maps. More precisely,
define $R_{A,\theta}:=\lambda_\theta(A)_*(A_\theta)_*(R)$ for $\theta\in V$. 
These currents belong to $\Cc (U)$. 
It is clear that $R_{A,\theta}$ depends continuously on $\theta\in
V$. We show that $(R_{A,\theta})$ defines a structural disc in
$\Cc(U)$, i.e. they are slices of some current $\Rc_A$. But first
observe that $R_{A,1}=R$, and $R_A:=R_{A,0}$ 
is independent
of $R$, and is defined by integration on the projective space $A(L)$. 
In other words, the structural disc $(R_{A,\theta})$, which contains
$R$, passes through a fixed
point $R_A$ when $R$ varies.

Define the meromorphic map $\Lambda_A:V\times\P^k\rightarrow \P^k$ by
$$\Lambda_A(\theta,x):=A_\theta^{-1} \lambda_\theta(A)^{-1}(x)$$
which is locally a holomorphic submersion outside $\{0\}\times
\P^k$. Then the positive closed $(p,p)$-current
$\Rc_A:=\Lambda_A^*(R)$ is well defined in $\Lambda_A^{-1}(U)$ which
is an open set in 
$(V\setminus\{0\})\times U$. Moreover, when $\theta\rightarrow 0$, the support
of $\Lambda_A^{-1}(U)$ clusters only on the set $\{0\}\times A(L)$. Then $\Rc_A$ is a
positive closed current on $(V\times U)\setminus  (\{0\}\times
A(L))$. Since the dimension of $\{0\}\times A(L)$ is strictly smaller than
the dimension of $\Rc_A$, the trivial extension of $\Rc_A$ across
$\{0\}\times A(L)$ is positive and closed \cite{HarveyPolking}. It is easy to check that the
slices $\langle \Rc_A,\pi_V,\theta\rangle$ are equal to $R_{A,\theta}$
for $\theta\not=0$. This family of currents is continuous and converges
to $R_{A,0}$ as $\theta\rightarrow 0$. We deduce that $\langle
\Rc_A,\pi_V,\theta\rangle = R_{A,\theta}$ for every $\theta$, see
Appendix \ref{appendix_slicing} below, and that
$(R_{A,\theta})$ defines a continuous structural disc in $\Cc(U)$.

We now introduce a smoothing. 
Let $\rho$ be a smooth positive probability measure 
with compact support in $\Aut(\P^k)$.
If $R$ is a current in $\P^k$ then 
the current $\int A_*(R)\d\rho(A)$, $A\in\Aut(\P^k)$, is smooth and converges to $R$ when $\rho$
tends to the Dirac mass at  $\id\in \Aut(\P^k)$. Fix $\rho$ with 
support in $W'$.
Define 
$$R_\theta:=\int R_{A,\theta}\d\rho(A).$$ 
It is clear that $(R_\theta)$ 
is also a continuous structural disc in $\Cc (U)$. They are slices of the current
$\Rc:=\int\Rc_A\d\rho(A)$. We also have $R_1=R$ and $R_0$  independent of $R$.

\begin{proposition} \label{prop_disc}
The current $R_\theta$ is smooth for every $\theta\in V\setminus\{1\}$. Moreover there exist constants $r>0$,
$c>0$ and $m>0$, independent of $R$, such that
\begin{enumerate}
\item $\|R_\theta-R_0\|_\infty\leq c|\theta|$ for $|\theta|<r$;
\item If $R$ is a continuous form and $M(R,\cdot)$ denotes the modulus of continuity of $R$, then
$\|R_\theta-R\|_\infty\leq c(|\theta-1|\|R\|_\infty+M(R,m|\theta-1|))$
for $|\theta-1|<r$. 
\end{enumerate}
\end{proposition}
\begin{proof} We have $R_\theta=\int A_*(A_\theta)_*R\d (\lambda_\theta)_*\rho(A)$. For $\theta\not=1$, 
since $(\lambda_\theta)_*\rho$ is a smooth function, $R_\theta$ is smooth.

1. Let $\Phi$ be a smooth $(k-p,k-p)$-form with compact support in
$U$. One needs to show that 
$|\langle R_\theta-R_0, \Phi\rangle|\leq c|\theta|\|\Phi\|_{\Lc^1}$. 
Using the identity $R_\theta=\int A_*(A_\theta)_*R\d
(\lambda_\theta)_*\rho(A)$, 
we get 
$$\langle R_\theta,\Phi\rangle = \langle R,\Phi'_\theta\rangle, \quad
{\rm where}\quad 
\Phi'_\theta :=   (A_\theta)^* \int A^*\Phi
\d(\lambda_\theta)_*\rho(A).$$ 
Hence, $\langle R_\theta-R_0,\Phi\rangle = \langle
R,\Phi'_\theta-\Phi'_0\rangle$ and
\begin{eqnarray*}
\Phi'_\theta -\Phi'_0 & :=  & (A_\theta)^* \int A^*\Phi
\d(\lambda_\theta)_*\rho(A) - (A_0)^* \int A^*\Phi \d\rho(A)\\
& = & (A_\theta)^* \int A^*\Phi \d\left((\lambda_\theta)_*\rho -\rho\right)(A)\\
& & +  (A_\theta)^* \int A^*\Phi \d\rho(A)- (A_0)^* \int A^*\Phi \d\rho(A).
\end{eqnarray*}
The form $\Phi'_\theta-\Phi'_0$ is defined and smooth in $\P^k\setminus I$. We only have to check that 
$\|\Phi'_\theta-\Phi'_0\|_{U,\infty}\leq c|\theta|\|\Phi\|_{\Lc^1}$,
i.e. to bound the last sum. 

Since 
$\|(\lambda_\theta)_*\rho-\rho\|_{\Cc ^1}
\leq\const|\theta|$, the first integral in the sum above defines a form with
$\Cc^1$-norm
bounded by  $\const|\theta|\|\Phi\|_{\Lc^1}$. It follows that
the first term of the sum is bounded 
by $\const|\theta|\|\Phi\|_{\Lc^1}$.
The last two integrals define a  form with $\Cc^1$-norm bounded by
$\const\|\Phi\|_{\Lc^1}$. We deduce that the rest of the sum is 
also bounded by $\const|\theta|\|\Phi\|_{\Lc^1}$ since 
$\|A_\theta-A_0\|_{\Cc ^1(U)}\leq \const|\theta|$.

2. There exists $m>0$ such that when $|\theta-1|<r$, $r$ small enough, and $A\in W'$ 
we have $\|\lambda_\theta(A)\circ A_\theta-\id\|_{\Cc ^1}\leq
m |\theta-1|$. Hence, 
$$\|R_{A,\theta}-R\|_\infty\leq \const 
(|\theta-1|\|R\|_\infty+M(R,m|\theta-1|)).$$ 
The estimate for $R_\theta$ follows from its definition.
\end{proof}


\section{Attracting current} \label{section_current}

In this section we  prove  Theorem \ref{th_attracting_current}.
Let $\Dc$ denote the set of all the currents $S$ obtained as limit values of $d^{-(k-p)n}(f^n)_*R_n$ 
with $R_n$ in $\Cc(U)$.
Observe that such a current $S$ has support in $\Ac$ and there exist $S_n\in \Dc$
such that $S=d^{-(k-p)n}(f^n)_*S_n$. Indeed, if $S=\lim d^{-(k-p)n_i}
(f^{n_i})_*R_{n_i}$ and if $S_n$ is a limit value of $d^{-(k-p)(n_i-n)}
  (f^{n_i-n})_*R_{n_i}$, $i\rightarrow\infty$, then $S=d^{-(k-p)n}(f^n)_*S_n$. In particular,
  $\Dc$ is a convex compact set.

\begin{proposition} \label{prop_convergence}
Let $\Phi$ be a real smooth $(k-p,k-p)$-form on $U$ such that $\ddc\Phi\geq 0$.
Then there exists a constant $c_\Phi$ such that $\langle S,\Phi\rangle\leq c_\Phi$
for every $S\in\Dc$ and 
if $R$ is a continuous form in   $\Cc (U')$,
then 
$$\lim_{n\rightarrow\infty}\langle d^{-(k-p)n}(f^n)_*R,\Phi\rangle = c_\Phi.$$
In particular, if $\ddc\Phi=0$ we have $\langle S,\Phi\rangle= c_\Phi$ for every 
$S\in\Dc$.
\end{proposition}
\begin{proof}
The case where $\Phi=\omega_\FS^{k-p}$ is trivial since $\langle
d^{-(k-p)n} (f^n)_* R_n,\Phi\rangle =1$. Replacing $\Phi$ by
$A\Phi+{1\over 2}\omega_\FS^{k-p}$, $A>0$, we can assume that $0\leq \Phi\leq \omega_\FS^{k-p}$ on $U'$. 
Then for $n\geq 1$ we have 
$\|d^{-(k-p)n} (f^n)^*\Phi\|_{U} \leq \|d^{-(k-p)n}
(f^n)^*\omega_\FS^{k-p}\|=1$
and $\ddc d^{-(k-p)n} (f^n)^*\Phi\geq 0$ on $U$. 

Let $(i_n)$ be a sequence of integers, $i_n>n$, and $R'_{i_n}$ be positive 
closed $(p,p)$-currents 
of mass 1 with 
support in $U$ such that $\langle d^{-(k-p)i_n} (f^{i_n})_* R'_{i_n},\Phi\rangle$ 
converges to a constant $c_\Phi$.
We choose $i_n$ and $R'_{i_n}$ so that $c_\Phi$ is the maximal value that we can 
obtain in this way. It is clear that $c_\Phi$ satisfies the inequality in the proposition. 
Since $0\leq \Phi\leq \omega_\FS^{k-p}$, we have $0\leq c_\Phi\leq 1$. 
Define $R_n:=d^{-(k-p)(i_n-n)}(f^{i_n-n})_*(R'_{i_n})$. Then
$\lim \langle d^{-(k-p)n}(f^n)_*R_n,\Phi\rangle =c_\Phi$. 
All the currents $R_n$ are supported in $U'$. 

Define the structural discs $(R_{n,\theta})$ associated to $R_n$ as in the 
last section. We have 
$$\varphi_n(\theta):=\langle d^{-(k-p)n} (f^n)_* R_{n,\theta},\Phi \rangle = 
\langle R_{n,\theta}, d^{-(k-p)n} (f^n)^*\Phi \rangle.$$
Since $R_{n,\theta}$ depends continuously on $\theta$ and since
$d^{-(k-p)n}(f^n)^*\Phi$ is smooth,
$\varphi_n(\theta)$ is a continuous function. Lemma
\ref{lemma_structure_coherent}, applied to $R_{n,\theta}$ and
to $d^{-(k-p)n}(f^n)^*\Phi$, 
implies that $\varphi_n(\theta)$ is
subharmonic on $\theta\in V$.  
On the other hand,   we have
$\limsup\varphi_n\leq c_\Phi$ (by definition of $c_\Phi$)
and
$$\lim_{n\rightarrow\infty} \varphi_n(1)=\lim_{n\rightarrow\infty} 
\langle d^{-(k-p)n} (f^n)_*R_n,\Phi\rangle = c_\Phi.$$
By maximum principle and Hartogs lemma, $\varphi_n$ converges in
$\Lc^1_\loc(V)$ to $c_\Phi$,
see \cite{Hormander}.
Proposition \ref{prop_disc} and the inequality $\|d^{-(k-p)n} (f^n)^*\Phi\|_{U}\leq 1$
imply that $|\varphi_n(\theta)-\varphi_n(0)|\leq c |\theta|$. 
We then deduce that $\lim\varphi_n(0)=c_\Phi$.

Consider now a continuous form $R\in \Cc (U')$. 
Define, as in the last section, the structural disc
$(R_\theta)$ associated to $R$ and continuous subharmonic functions
$$\psi_n(\theta):=\langle d^{-(k-p)n} (f^n)_* R_{\theta},\Phi \rangle = 
\langle R_{\theta}, d^{-(k-p)n} (f^n)^*\Phi
\rangle.$$
Since $R_{0}=R_{n,0}$ we have $\psi_n(0)=\varphi_n(0)$. 
Hence, $\lim \psi_n(0)=c_\Phi$. 
This, combined with the inequality $\limsup\psi_n\leq c_\Phi$, implies 
that $\psi_n\rightarrow c_\Phi$ in $\Lc^1_\loc(V)$. On the other hand, 
Proposition \ref{prop_disc} gives that 
$\lim_{\theta\rightarrow 1}(\sup_n |\psi_n(\theta)-\psi_n(1)|)=0$. 
Therefore, $\lim \psi_n(1)=c_\Phi$. 
Since $R_{1}=R$, we obtain $\lim \langle d^{-(k-p)n}(f^n)_*R,\Phi\rangle=c_\Phi$.

When $\ddc\Phi=0$, the inequality
in the proposition applied to $\pm \Phi$ yields $\langle S,\Phi\rangle =c_\Phi$
 for every $S\in\Dc$.
\end{proof}

The last proposition and Lemma \ref{lemma_density} 
show that if $R$ is a continuous positive closed $(p,p)$-form of mass 1
with support in $U$ then $d^{-(k-p)n}(f^n)_*R$ converges to a current $\tau$
which is independent
of $R$ (we can choose $U'$ so that $\supp(R)\subset U'$). This current is given by
$$\langle\tau,\Phi\rangle := c_\Phi.$$
This is the convergence in Theorem \ref{th_attracting_current}

\begin{remark} \rm
If $(R_n)\subset \Cc(U)$ such that $R_n\leq cR$, $c\geq 0$, then 
$d^{-(k-p)n} (f^n)_* R_n$ converges to $\tau$. Indeed, Proposition
\ref{prop_convergence}
implies that 
$$\limsup \langle d^{-(k-p)n} (f^n)_* R_n,\Phi\rangle \leq c_\Phi$$
and 
$$ \limsup \langle d^{-(k-p)n} (f^n)_* (cR-R_n),\Phi\rangle \leq (c-1)c_\Phi.$$ 
Consider the sum of these inequalities. 
Since  $\lim \langle d^{-(k-p)n} (f^n)_* R,\Phi\rangle =c_\Phi$, we deduce that
$\lim \langle d^{-(k-p)n} (f^n)_* R_n,\Phi\rangle =c_\Phi$.
\end{remark}

The following proposition shows that one needs to test only one form in order to 
check the convergence toward $\tau$.

\begin{proposition} \label{prop_one_test} 
Let $\Phi$ be a real smooth $(k-p,k-p)$-form on $U$ such that $\ddc\Phi$ is strictly positive.
Let $(R_{n_i})$ be a sequence of currents in $\Cc (U)$. 
Then $d^{-n_i}(f^{n_i})_*R_{n_i}$ converges to $\tau$ if and
only if $\lim \langle d^{-n_i}(f^{n_i})_*R_{n_i},\Phi\rangle=\langle\tau,\Phi\rangle$.
\end{proposition}
\begin{proof} 
Let $\Psi$ be a real smooth $(k-p,k-p)$-form with compact support in $U$. Let $c>0$ such that $\ddc(c\Phi\pm\Psi)
\geq 0$. Proposition \ref{prop_convergence} implies that 
$$\limsup_{n\rightarrow\infty} \langle d^{-n_i}(f^{n_i})_*R_{n_i},c\Phi\pm\Psi\rangle
\leq\langle\tau,c\Phi\pm\Psi\rangle.$$
If $\lim \langle d^{-n_i}(f^{n_i})_*R_{n_i},\Phi\rangle=\langle\tau,\Phi\rangle$ then
$\lim \langle d^{-n_i}(f^{n_i})_*R_{n_i},\Psi\rangle=\langle\tau,\Psi\rangle$ for every $\Psi$.
The proposition follows.
\end{proof}

We have the following Corollary.

\begin{corollary} \label{cor_laminarity}
Let $L'$ be a generic projective space of dimension $k-p$, close enough to $L$. Then
there exists an increasing sequence $(n_i)$ of integers such that
$d^{-n_i}(f^{n_i})_*[L']\rightarrow\tau$. In particular, $\tau$ is
weakly laminar.
\end{corollary}
\begin{proof} 
Let $\rho$ be a smooth function as above such that $\rho>0$ on a
neighbourhood of $\id\in\Aut(\P^k)$. We want to prove the convergence
for $L':=A(L)$ and for $\rho$-almost every $A\in\Aut(\P^k)$.
Let $\Phi$ be as in Proposition \ref{prop_one_test}. 
We have seen that 
$$\lim_{n\rightarrow\infty} \langle d^{-(k-p)n}(f^n)_*R_0,
\Phi\rangle = c_\Phi.$$ 
Since $R_0=\int[A(L)] \d \rho(A)$, we have
$$\lim_{n\rightarrow\infty} \int \langle d^{-(k-p)n} (f^n)_*[A(L)],
\Phi\rangle \d\rho(A) =c_\Phi.$$
On the other hand, the integrals $\langle d^{-(k-p)n} (f^n)_*[A(L)],
\Phi\rangle$ are unifomly bounded with respect to $n$ and $A$, and 
$$\limsup_{n\rightarrow\infty}   \langle d^{-(k-p)n} (f^n)_*[A(L)],
\Phi\rangle\leq c_\Phi.$$ 
Therefore, 
for $\rho$-almost every $A$ there exists $(n_i)$ 
such that 
$$\lim_{i\rightarrow\infty}\langle d^{-n_i} (f^{n_i})_*[L'],\Phi\rangle
=c_\Phi, \quad \mbox{where } L':=A(L).$$ 
This, together with Proposition \ref{prop_one_test}, implies that  
$d^{-n_i}(f^{n_i})_*[L']\rightarrow\tau$.

Let $L_i$ be holomorphic images of $\P^{k-p}$ in $\P^k$ and let
$c_i$ be positive numbers. Assume that the sequence of currents
$c_i[L_i]$ converges. Theorem 5.1 in \cite{Dinh} implies that the
limit is a weakly laminar current. Then $\tau$ is weakly laminar.  
\end{proof}

The following proposition completes the proof of Theorem \ref{th_attracting_current}.

\begin{proposition} \label{prop_extremality}
The current $\tau$ is extremal in $\Dc$. In particular, it is extremal
in the cone of invariant positive closed $(p,p)$-currents supported on $\Ac$.
\end{proposition}
\proof
Recall that $\Dc$ is a convex compact set.
Assume that $\tau=\lambda S_1+(1-\lambda) S_2$ with $S_1$ and $S_2$ in
$\Dc$. Then Proposition \ref{prop_convergence} implies that
$\langle S_1,\Phi\rangle \leq c_\Phi$ and $\langle S_2,\Phi\rangle
\leq c_\Phi$.
We also have
$$c_\Phi= \langle \tau,\Phi\rangle = \lambda\langle S_1,\Phi\rangle
+ (1-\lambda) \langle S_2,\Phi\rangle.$$
Hence $\langle S_1,\Phi\rangle =\langle S_2,\Phi\rangle =c_\Phi$. It
follows that $S_1=S_2=\tau$. 
\endproof

We have the following version of Theorem \ref{th_attracting_current}
where the form $R'$ is continuous but not necessarily positive and
closed.

\begin{theorem} \label{th_current_bis}
Let $S$ be a positive closed $(p,p)$-current in $\Cc(U)$. Let $(n_i)$ be an
increasing sequence of integers such that $d^{-(k-p)n_i}(f^{n_i})_*S$
converges to $\tau$. If $\varphi$ is a continuous function on $U$, then 
 $d^{-(k-p)n_i}(f^{n_i})_*(\varphi S)$ converges weakly to $c_{\varphi S} \tau$, where 
$c_{\varphi S}:=\langle S\wedge T^{k-p},\varphi\rangle$.
In particular,  $d^{-(k-p)n}(f^n)_*(\varphi \tau)$ converges weakly to
$c_{\varphi\tau} \tau$. 
If $R$ is as in Theorem \ref{th_attracting_current} and  $R'$ is 
a continuous real $(p,p)$-form, not necessarily positive and closed, 
such that $-cR\leq R'\leq cR$, $c>0$,
then $d^{-(k-p)n}(f^n)_*R'$ converges weakly to $c_{R'} \tau$, where 
$c_{R'}:=\langle T^{k-p},R'\rangle$.
\end{theorem} 

We first prove the following proposition which can be easily extended 
to the case of meromorphic maps or correspondences, see also
\cite{BedfordSmillie, Sibony} and
\cite[Remarque 5.9]{Dinh}.

\begin{proposition} \label{prop_d_convergence}
Let $\varphi$ be a real-valued smooth function and let $S$ be a positive
closed $(m,m)$-current on $\P^k$, $0\leq m\leq k-1$. Then 
$$\|\ddc \big(d^{-(k-m)n} (f^n)_*(\varphi S)\big)\|= O(d^{-n}) \quad \mbox{and} \quad 
\|\d \big(d^{-(k-m)n}(f^n)_*(\varphi S)\big)\|= O(d^{-n/2}).$$
In particular, if $R'$ is a smooth $(m,m)$-form on $\P^k$ then
$$\|\ddc \big(d^{-(k-m)n} (f^n)_*R'\big)\|= O(d^{-n}) \quad \mbox{and} \quad 
\|\d \big(d^{-(k-m)n}(f^n)_*R'\big)\|= O(d^{-n/2}).$$
\end{proposition}
\proof
Observe that the first estimate is classical. We can write $\ddc\varphi$ as a  linear combination of
positive closed $(1,1)$-forms. Therefore,
$\ddc\varphi \wedge S$ can be written as a linear combination of positive closed
$(m+1,m+1)$-currents. 
On the other hand, we have seen in the Introduction that if $S'$ is a positive closed $(m+1,m+1)$-current,
then $\|(f^n)_*S'\|=
d^{(k-m-1)n}\|S'\|$. Hence
$$\|\ddc \big(d^{-(k-m)n} (f^n)_*(\varphi S)\big)\|=
\|\big(d^{-(k-m)n} (f^n)_*(\ddc \varphi \wedge S)\big)\|=O(d^{-n}).$$

For the second estimate in the proposition, we use,  as
in \cite{Dinh, DinhSibony5}, the product map
$F:\P^k\times\P^k\rightarrow \P^k\times\P^k$ with
$F(x,y)=(f(x),f(y))$. 
A simple computation on cohomology implies that $\|(F^n)_*\widetilde
S\|\lesssim d^{(2k-s)n}$ for every positive closed $(s,s)$-current
$\widetilde S$ in $\P^k\times\P^k$.
We will show that 
$$\|\partial \big(d^{-(k-m)n}(f^n)_*(\varphi S)\big)\|= O(d^{-n/2}).$$
We also need an analogous inequality with $\dbar$ instead of
$\partial$, which  is
proved in the same way.
Let $\Phi$ be a smooth $(k-m-1,k-m)$-form on $\P^k$ such that
$\|\Phi\|_\infty=1$. It is sufficient to check that 
$$|\langle (f^n)_*(\partial\varphi \wedge S),\Phi\rangle|
=|\langle \partial (f^n)_*(\varphi S),\Phi\rangle| \leq c
d^{(k-m-{1\over 2})n},$$
where $c>0$ is independent of $\Phi$. 

Consider the following forms and positive closed current in $\P^k\times\P^k$
$$\widehat\phi:=i\partial\varphi\otimes\dbar\varphi,\quad \widehat S =
S\otimes S\quad \mbox{and} \quad \widehat\Phi=\Phi\otimes\overline
\Phi.$$ 
Then the desired inequality is equivalent to
$$|\langle (F^n)_*(\widehat\phi\wedge \widehat
S),\widehat\Phi\rangle|\leq c^2 d^{(2k-2m-1)n}.$$
Since $\widehat\phi$ is smooth, it can be written as a linear
combination of positive $(1,1)$-forms which are bounded by closed ones. Then $\widehat\phi\wedge
\widehat S$ can be bounded by a positive closed
$(2m+1,2m+1)$-current. The last inequality follows from the previous
observation on the mass of $(F^n)_*(\widetilde S)$.

When $R'$ is smooth, it can be written as a linear combination
of forms of type $\varphi R$, where $\varphi$ is a smooth function and
$R$ is a positive closed $(m,m)$-form. The rest of the
proposition follows.

Finally, note that in the proposition, one can replace $\ddc\phi$,
$\d\phi$, $\ddc R$, $\d R$ by continuous forms not necessarily closed. 
\endproof

\noindent
{\bf Proof of Theorem \ref{th_current_bis}.}
Replacing $\varphi$ by $\varphi+A$, $A>0$, we can assume that
$\varphi$ is strictly positive. Since $\varphi$ can be uniformly
approximated by smooth functions, we can also assume that $\varphi$ is
smooth.

Let $\tau'$ be the limit
of $d^{-(k-p)m_i}(f^{m_i})_*(\varphi S)$ for  a
subsequence $(m_i)$  of $(n_i)$. Observe that $\tau'$ is positive.
By Proposition
\ref{prop_d_convergence}, $\tau'$ is closed.
We also have
\begin{eqnarray*}
\|\tau'\| & = & \int \tau'\wedge \omega_\FS^{k-p}=\lim_{i\rightarrow\infty}
\int d^{-(k-p)m_i}(f^{m_i})_*(\varphi S) \wedge \omega_\FS^{k-p}  \\
& = & \lim_{i\rightarrow\infty}
\int \varphi S\wedge d^{-(k-p)m_i}(f^{m_i})^*\omega_\FS^{k-p}= \int
\varphi S\wedge
T^{k-p}=c_{\varphi S},
\end{eqnarray*}
see \cite{Sibony} for the last limit.
Moreover, if $\tau'_n$ is a limit
value of $d^{-(k-p)(m_i-n)}(f^{m_i-n})_*(\varphi S)$ then $\tau_n'$ is
closed and $\tau'=d^{-(k-p)n}(f^n)_*\tau'_n$. 
It follows that $c_{\varphi S}^{-1} \tau'$ is an element of $\Dc$. 
On the other hand, we have $\tau'\leq \lim d^{-(k-p)n_i}
(f^{n_i})_*(cS)=c\tau$ for $c>0$ such that $\varphi\leq c$. 
Proposition \ref{prop_extremality} implies
that $\tau'=c_{\varphi S}\tau$. Therefore, $d^{-(k-p)n_i}(f^{n_i})_*(\varphi
S)$ converges to $c_{\varphi S} \tau$.

The assertion on $\varphi\tau$ is clear.

Arguing as in Lemma \ref{lemma_density}, we
can regularize $R$ and obtain 
 a positive closed form $\widetilde R$ of mass 1 with compact support
 in $U$, strictly positive in a neighbourhood
$U''$ of $\supp(R')$. Replace $R$ by $\widetilde R$. Then, 
since $R'$ can be uniformly
approximated by smooth forms in $U''$, we can assume that $R'$ is
smooth. We always have $-cR\leq R'\leq c R$ for some constant
$c>0$. Now,  by Theorem
\ref{th_attracting_current}, we can replace $R'$ by $R'+cR$ and assume that $R'$ is positive.
Then the previous proof for $\varphi S$ is valid for $R'$ in the
present context. We use
here Theorem \ref{th_attracting_current} and the second part of Proposition
\ref{prop_d_convergence}. 
\hfill $\square$


\section{Equilibrium measure} \label{section_measure}

In this section we prove Theorem \ref{th_measure}. We have seen in Corollary
\ref{cor_entropy} that the entropy of $\nu$ satisfies $h(\nu)\leq (k-p)\log d$. In order to prove the
inequality $h(\nu)\geq (k-p)\log d$, we can follow the proof of
Bedford-Smillie for H{\'e}non maps \cite{BedfordSmillie}, see also
\cite{Yomdin, deThelin2}. For this purpose one only needs the following lemma.

\begin{lemma} Let $L'$ be a generic projective space of dimension $k-p$, close enough to $L$. Then
there exists an increasing sequence $(n_i)$ of integers such that
$${1\over n_i} \sum_{j=0}^{n_i-1}d^{-(k-p)n_i}(f^j)_*[L']\wedge
(f^{n_i-j})^*\omega_\FS^{k-p}\rightarrow\nu.$$ 
\end{lemma}
\proof
As in Corollary \ref{cor_laminarity}, we prove that if $L'$ is generic, 
$${1\over n_i} \sum_{j=0}^{n_i-1}d^{-(k-p)j}(f^j)_*[L']\rightarrow\tau.$$ 
Then
$${1\over n_i} \sum_{j=0}^{n_i-1}d^{-(k-p)j}(f^j)_*[L']\wedge
T^{k-p}\rightarrow\nu,$$
see  \cite{Demailly,
  FornaessSibony3, Sibony}. One only needs the following fact, applied
to $S=d^{-(k-p)j}(f^j)_*[L']$ in order
to complete the proof.

Write $d^{-n} (f^n)^*\omega_\FS=\omega_\FS+\ddc g_n$, where $g_n$ is a smooth function 
converging uniformly to the Green function $g$ of $f$ as $n\rightarrow\infty$, see
e.g. \cite{Sibony}. We have $T=\omega_\FS+\ddc g$. If $S$ is a
positive closed current of mass 1 and if $\Phi$ is a fixed smooth test
form of right degree, it is easy to show that 
$$\langle S\wedge  d^{-(k-p)n} (f^n)^*\omega_\FS^{k-p},\Phi\rangle -
\langle S\wedge T^{k-p},\Phi\rangle$$ 
converges to 0 uniformly on $S$, see \cite{Demailly,
  FornaessSibony3, Sibony}.
\endproof

Now, we have to prove the mixing. We follows the ideas of
\cite{BedfordSmillie, Sibony}, see \cite{DinhSibony5} for another approach.
Let $\varphi$ and $\psi$ be real-valued smooth functions. We want to check that 
$$\langle \nu, \varphi(\psi\circ f^n)\rangle 
\rightarrow \langle \nu, \varphi\rangle\langle\nu,\psi\rangle.$$
We have
$$\langle \nu,\varphi(\psi\circ f^n)\rangle = 
\langle \varphi \tau\wedge T^{k-p},\psi\circ f^n\rangle = 
\langle (f^n)_*(\varphi \tau\wedge T^{k-p}),\psi\rangle.$$
Define $c_\varphi:=\langle\nu,\varphi\rangle$. 
Then the following proposition, applied to $m=k-p$,  implies the result.

\begin{proposition}
For all $0\leq m\leq k-p$, the current
$d^{-(k-p-m)n}(f^n)_*(\varphi\tau\wedge T^m)$ converges weakly to $c_\varphi
  \tau\wedge T^m$ as $n\rightarrow\infty$. 
\end{proposition}
\proof
The proof uses an induction on $m$. For $m=0$ the proposition is
deduced from Theorem \ref{th_current_bis}. Assume the proposition
for $m$, $0\leq m\leq k-p-1$. Let $\Psi$ be a smooth test form
of right bidegree. Define $S':=\varphi\tau\wedge T^m$.
We have
\begin{eqnarray*}
\lefteqn{d^{-(k-p-m-1)n}\langle (f^n)_*(\varphi\tau\wedge T^{m+1}),
 \Psi\rangle}\\
 & = & d^{-(k-p-m-1)n}\langle S'\wedge T, (f^n)^*\Psi\rangle \\
&=& \lim_{s\rightarrow\infty} d^{-(k-p-m-1)n}\langle
S'\wedge d^{-n-s}(f^{n+s})^*\omega_\FS,(f^n)^*\Psi\rangle \\
&=& \lim_{s\rightarrow\infty} d^{-(k-p-m-1)n}\langle
S',d^{-n-s}(f^{n+s})^*\omega_\FS\wedge (f^n)^*\Psi\rangle \\
&=& \lim_{s\rightarrow\infty} d^{-(k-p-m)n}\langle
(f^n)_*S',d^{-s}(f^s)^*\omega_\FS\wedge \Psi\rangle \\
&=& \lim_{s\rightarrow\infty} d^{-(k-p-m)n}\langle
(f^n)_*S',(\omega_\FS+\ddc g_s)\wedge \Psi\rangle 
\end{eqnarray*}
The last integral is equal to 
\begin{eqnarray*}
d^{-(k-p-m)n}\lefteqn{\langle  (f^n)_*S',\Psi\wedge \omega_\FS +g_s\ddc\Psi \rangle
 + d^{-(k-p-m)n}\langle \ddc [ (f^n)_*S'], g_s\Psi \rangle} \\ 
&& +  d^{-(k-p-m)n}\langle \d [ (f^n)_*S'], 
g_s\dc \Psi \rangle - d^{-(k-p-m)n} \langle \dc [ (f^n)_*S'], g_s\d \Psi \rangle.
\end{eqnarray*}
Since the $g_s$ are uniformly bounded, Proposition \ref{prop_d_convergence} implies that the  
last three terms tend to 0 when $n$ and $s$ tend to infinity. Since $g_s$
converges uniformly to $g$, the induction hypothesis implies that the
first term converges to
$$c_\varphi\langle \tau\wedge T^m, \Psi \wedge \omega_\FS+ g\ddc\Psi \rangle = c_\varphi
\langle \omega_\FS\wedge\tau \wedge T^m+ \ddc (g\tau\wedge T^m),
\Psi\rangle=c_\varphi\langle \tau\wedge T^{m+1},\Phi\rangle$$
since $\omega_\FS+\ddc g =T$. This completes the proof.
\endproof

\section{Remarks and questions} \label{section_remark}

We consider the case where $p=1$ and $\P^k\setminus U$ is convex in
$\P^k\setminus L\simeq \C^k$. The fact that $\tau$ is a current of
bidegree $(1,1)$ allows us to prove some further properties of the
attracting set $\Ac$. 

The following proposition shows that for a generic point in $\P^k$ 
most of its preimages are out of $U$. 

\begin{proposition} Let $\Leb$ denote the Lebesgue measure in $\P^k$.
Then $\|(f^n)^*\Leb\|_{U}=o(d^n)$.
In particular, for every $\epsilon>0$ 
$$\lim_{n\rightarrow \infty}\!\ \Leb\{a\in \P^k, \ \# f^{-n}(a)\cap U\geq \epsilon d^n\}
=0.$$ 
\end{proposition}
\begin{proof} 
Observe that $U$ admits a neighbourhood $U''$ satisfying the same
properties, 
i.e. $f(U'')\Subset U''$ and $\P^k\setminus U''$ is
convex. We replace $U$ by $U''$. Then we can replace 
$U$ in the proposition by a compact set $K\subset U$.

Let $S$ be a smooth $(k-1,k-1)$-form on $\P^k\setminus I$ such that
$\ddc S=\Leb$. Here, we identify $\Leb$ with a volume form. Let 
$\varphi$ be 
a positive smooth function with support in $U$ such that $\varphi=1$ on $K$.

Observe that since $\P^k\setminus U$ is convex, we can find a smooth
form $R\in\Cc(U)$ strictly positive on the support of $\ddc\varphi$. 
Theorem \ref{th_current_bis} applied to $R':=\ddc\varphi$ and to $R$ yields
$d^{-n}(f^n)_*R'\rightarrow 0$. Hence
$$\langle d^{-n}(f^n)^*\Leb,\varphi\rangle = 
\langle d^{-n}\ddc (f^n)^*S,\varphi\rangle = 
\langle S, d^{-n}(f^n)_*R'\rangle
\rightarrow 0.$$
The first assertion follows.
 
Now observe that $\#f^{-n}(a)\cap K\geq \epsilon d^n$ 
if and only if $(f^n)_*{\bf 1}_{K}(a)\geq \epsilon d^n$. 
Hence
\begin{eqnarray*}
\lefteqn{\epsilon d^n \Leb\{a\in \P^k, \ \# f^{-n}(a)\cap K\geq
  \epsilon d^n\}}\\ 
& \leq &  
\langle \Leb, (f^n)_*{\bf 1}_{K}\rangle
 =  \langle (f^n)^*\Leb,{\bf 1}_{K}\rangle  =  \|(f^n)^*\Leb\|_{K}.
\end{eqnarray*}
We conclude using the first assertion in the proposition.
\end{proof}

Let $R$ be a positive closed $(1,1)$-current of mass 1 in $U$. 
We call {\it canonical quasi-potential}
of $R$ the unique quasi-psh function $g_R$ on $\P^k$ such that $g_R(I)=0$ and 
$\ddc g_R = R-\omega_\FS$. 
Denote by $g_\tau$ the canonical quasi-potential of $\tau$.
The following result gives a refined version of Theorem \ref{th_attracting_current}
for the bidegree $(1,1)$ case.

\begin{theorem} \label{th_canonical_potential}
Let $S$ be an element of $\Dc$ and $g_S$ its canonical quasi-potential.
Then $g_S\leq g_\tau$ on $\P^k$ and $g_S=g_\tau$ on $U$. 
If the jacobian of $f$ satisfies $|\Jac(f)|<1$ on $U$ then $\tau$ is the unique invariant positive
closed $(1,1)$-current with mass $1$ supported in $U$. 
\end{theorem}
\begin{proof}
We prove that $g_S=g_\tau$ in $\P^k\setminus \overline U$. 
Consider a $(k-1,k)$-form $\Psi$ with support in $\P^k\setminus\overline U$.
We first show that $\langle g_S-g_\tau,\partial\Psi\rangle=0$.  
By Hodge theory, $\partial\Psi$ is $\ddc$-exact in $\P^k$. 
Write $\Psi=\ddc\Phi$. We have
$\langle g_S-g_\tau,\partial\Psi \rangle=\langle S-\tau,\Phi\rangle$. Since 
$\ddc\Phi=0$ on $U$, Proposition \ref{prop_convergence} 
implies that the last integral vanishes. Hence, $\partial(g_S-g_\tau)=0$ 
on $\P^k\setminus \overline
U$. Since $g_S-g_\tau$ is real-valued, $g_S-g_\tau$ is constant in $\P^k\setminus \overline U$.
This constant is zero since $g_S(I)=g_\tau(I)=0$.
Since  $S$ and $\tau$ have supports in $\Ac$, 
$g_S$ and $g_\tau$ are pluriharmonic on  $\P^k\setminus \Ac$. 
We deduce that $g_S=g_\tau$ on the connected component $W$ of
$\P^k\setminus\Ac$ which contains $\P^k\setminus U$. 

Let $\Omega$ be a smooth positive $(k,k)$-form with support in a small
neighbourhood of $\Ac$. 
We show that
$\langle g_S-g_\tau,\Omega\rangle$ is negative. 
There
exists a smooth real $(k-1,k-1)$-form $\Phi$ on $U$ such that $\ddc\Phi=\Omega$. 
Since $g_S-g_\tau=0$ on $W$,  we have $\langle g_S-g_\tau,\Omega\rangle = 
\langle g_S-g_\tau,\ddc\Phi\rangle
= \langle S-\tau,\Phi\rangle$. Since $\ddc\Phi$ is positive, 
Proposition \ref{prop_convergence} implies
that the last integral is negative. Hence, $g_S\leq g_\tau$.

Define
$\Ac':=\P^k\setminus W$. We show that
$f(\Ac')\subset \Ac'$.
We have
$f(\Ac')\Subset U$ and $\partial \Ac'\subset \Ac\subset\Ac'$. 
Since $f$ is open, $\partial f(\Ac')\subset f(\partial \Ac')\subset f(\Ac)\subset \Ac$ and 
then $f(\Ac')\cap W=\varnothing$.
We deduce from the definition of $W$ that $f(\Ac')\subset\Ac'$.
If the jacobian of $f$ satisfies $|\Jac(f)|<1$ on $U$ then  $\Ac'$ has
Lebesgue measure zero and $g_S=g_\tau$ almost
everywhere. Therefore $g_S=g_\tau$ everywhere since they are quasi-psh, and $\Dc$ is reduced
to $\{\tau\}$.
\end{proof}

\begin{remark}  \rm 
The hypothesis in Theorem \ref{th_canonical_potential} is stable under small
pertubations on $f$. In particular, small pertubations of polynomial
maps give us examples satisfying this condition. For such a map, since
the sum of Lyapounov exponents associated to $\nu$ is equal to
$\langle \nu, \log |\Jac(f)|\rangle$ which is negative, $\nu$ admits a
strictly negative Lyapounov exponent. In the dimension 2 case, de Th{\'e}lin
proved that $\nu$ has a strictly positive exponent \cite{deThelin2}; then $\nu$ is
hyperbolic.
\end{remark}

The following questions are still open for the general case.

\medskip

\noindent
{\bf Questions (with de Th{\'e}lin).} {\bf 1.} 
Is $\tau$ the unique (invariant) positive closed $(p,p)$-current of mass 1 with support in $\Ac$ ?
Is $\Ac$ always of Lebesgue measure zero ? Is the quasi-potential of
$\tau$ unbounded
in the bidegree $(1,1)$ case ? 

{\bf 2.} Is $\nu$ the unique measure of maximal entropy on $\Ac$ ? 
Is it exponentially mixing and singular with respect to the Lebesgue measure ?

{\bf 3.} 
Is $T$ laminar in $U$ ? Does $\nu$ have a product structure ? 
Are periodic saddle points equidistributed on the support of $\nu$ ?

{\bf 4.} Does $\nu$ have strictly negative Lyapounov exponents ? 
Even when we assume that, the previous questions are open.


\begin{appendix}

\section{Slicing of horizontal currents} \label{appendix_slicing}

In this appendix we give some results on the slicing of currents 
used in this paper.

Consider the open set  $U\subset \P^k$ as in the Introduction, the
manifold $V$ of dimension $m$ and the canonical projections
$\pi_V:V\times U\rightarrow V$, $\pi_U:V\times U\rightarrow U$ as in 
Section \ref{section_structure}.  Let $\Rc$ be a 
positive closed  $(p,p)$-current on $V\times U$.
The slicing theory
\cite{Federer, HarveyShiffman} allows us to define slices $\langle
\Rc,\pi_V,\theta\rangle$ for almost every $\theta\in V$, see also
\cite{DinhSibony5}. 
We can consider  $\langle
\Rc,\pi_V,\theta\rangle$  as the intersections of $\Rc$ with the
currents of integration on $\pi_V^{-1}(\theta)$. They are positive
closed $(p,p)$-currents on $\{\theta\}\times
U$. We often identify them to $(p,p)$-currents in $U$ or
consider them as currents on $V\times U$, of bidegree $(p+m,p+m)$.  

Slicing is the generalization of restriction of forms to level sets of
holomorphic maps.
When $\Rc$ is a continuous form, 
$\langle \Rc, \pi_V,\theta\rangle$ is simply the restriction of $\Rc$ to $\pi_V^{-1}(\theta)$. 
When $\Rc$ is the current of integration on an analytic subset $X$ of
$V\times U$, 
$\langle \Rc,\pi_V,\theta\rangle$ is the current of integration on the analytic set 
$X\cap \pi_V^{-1}(\theta)$ for $\theta$ generic. 
If $\Omega$ is a smooth form
of maximal degree with compact support in $V$ and if $\Psi$
is a smooth $(k-p,k-p)$-form in $V\times U$, then
$$\int_V \langle \Rc,\pi_V,\theta\rangle (\Psi)\Omega(\theta)=
\langle \Rc\wedge \pi_V^*(\Omega),\Psi\rangle.$$ 
Here, we need that $\pi_V$ is proper on $\supp(\Psi)\cap\supp(\Rc)$.
In particular, this holds when $\Psi$ has
compact support. 

Let $y$ denote the coordinates in a chart of $V$ and $\lambda_V$ the
associated standard volume form.
Let $\psi(y)$ be a positive smooth function with compact support
such that $\int\psi\lambda_V=1$. Define
$\psi_\epsilon(y):=\epsilon^{-2m}\psi(\epsilon^{-1} y)$ and 
$\psi_{\theta,\epsilon}(y):=\psi_\epsilon(y-\theta)$ (the measures $\psi_{\theta,\epsilon}\lambda_V$
approximate the Dirac mass at $\theta$). Then, for every smooth test $(k-p,k-p)$-form $\Psi$
with compact support in $V\times U$ 
one has
$$\langle \Rc,\pi_V,\theta\rangle (\Psi)=\lim_{\epsilon\rightarrow 0}
\langle \Rc\wedge \pi_V^*(\psi_{\theta,\epsilon}\lambda_V),\Psi\rangle$$
when $\langle \Rc,\pi_V,\theta\rangle$ exists.
This property holds for all choice of  $\psi$ and also for $\Psi$ such that $\pi_V$
is proper on $\supp(\Psi)\cap\supp(\Rc)$.
Conversely, when the previous limit exists and is independent of  $\psi$, 
it defines $\langle \Rc,\pi_V,\theta\rangle$ and one says 
that $\langle \Rc,\pi_V,\theta\rangle$ {\it is well defined}.

Our main result of this section is the following proposition.

\begin{proposition} \label{prop_slice}
Assume that $\Rc$ is horizontal, i.e. $\pi_U(\supp(\Rc))\Subset U$.
Then the slice $\langle
\Rc,\pi_V,\theta\rangle$ exists for every $\theta\in V$. Moreover, if
$\Psi$ is a real smooth $(k-p,k-p)$-form on $V\times U$ such that
$\ddc \Psi\geq 0$ then $\langle\Rc,\pi_V,\theta\rangle(\Psi)$ defines
a psh function on $\theta\in V$. If $\ddc \Psi= 0$ 
then $\langle\Rc,\pi_V,\theta\rangle(\Psi)$ is pluriharmonic.
\end{proposition}
\proof
Since the problem is local, we can assume that $V$ is a ball
in $\C^m$ and use the above notation.  Let $\Psi$ be a real smooth 
$(k-p,k-p)$-form on $V\times U$. We want to prove that 
$$\lim_{\epsilon\rightarrow 0}
\langle \Rc\wedge
\pi_V^*(\psi_{\theta,\epsilon}\lambda_V),\Psi\rangle$$
exists and does not depend on the choice of $\psi$.
An analogous result as Lemma \ref{lemma_density} allows us to assume
that $\ddc\Psi\geq 0$. 
Since $\pi_V$ is proper on the support
of $\Rc$, the current $\varphi:=(\pi_V)_*(\Rc\wedge \Psi)$ is well
defined and has bidegree $(0,0)$. We also have
$\ddc\varphi=(\pi_V)_*(\Rc\wedge \ddc \Psi) \geq
0$. Then $\varphi$ is (equal to) a psh function on $V$. We have
$$\lim_{\epsilon\rightarrow 0}
\langle \Rc\wedge
\pi_V^*(\psi_{\theta,\epsilon}\lambda_V),\Psi\rangle =
\lim_{\epsilon\rightarrow 0}
\int_V \varphi \psi_{\theta,\epsilon}\lambda_V.$$
A classical property of psh functions implies that the last limit is
equal to $\varphi(\theta)$. The proposition follows.  
\endproof

\begin{remark} \label{rk_slice}
\rm
If $\Phi$ is a real smooth $(k-p,k-p)$-form on $U$ such that
$\ddc\Phi\geq 0$, then
$\langle\Rc,\pi_V,\theta\rangle(\Phi)$
is also psh on $V$. Indeed, we can apply Proposition \ref{prop_slice}
to $\Psi:=\pi_U^*(\Phi)$.
\end{remark}

\begin{proposition} \label{prop_slice_mass}
Under the hypothesis of Proposition \ref{prop_slice}, the mass of
$\langle \Rc,\pi_V,\theta\rangle$ is independent of $\theta$.
\end{proposition}
\proof
Consider the positive closed form $\Phi:=\omega_\FS^{k-p}$. Define as
above $\Psi:=\pi_U^*(\Phi)$ and $\varphi:=(\pi_V)_*(\Rc\wedge
\Psi)$. Then $\varphi$ is closed. It follows that $\varphi$ is a
constant function. Hence the mass of   $\langle
\Rc,\pi_V,\theta\rangle$, 
which is equal to $\langle
\Rc,\pi_V,\theta\rangle(\Phi)=\varphi(\theta)$, does not depend on $\theta$.
\endproof

\end{appendix}

 \noindent
{\bf Acknowledgments.} I would like to thank Henry de Th{\'e}lin,
Viet-Anh Nguyen and Nessim Sibony for their remarks.

Tien-Cuong Dinh \\
Institut de Math{\'e}matique de Jussieu \\
Plateau 7D, Analyse Complexe \\
175 rue du Chevaleret \\
75013 Paris, France \\
{\tt dinh@math.jussieu.fr}
\\

\end{document}